\documentclass[12pt,a4paper,oneside]{article}
\usepackage{amsmath,amstext,amssymb,amsopn,amsthm}

\theoremstyle{plain}
\newtheorem{thm}{Theorem}[section]

\newtheorem{lem}[thm]{Lemma}

\theoremstyle{definition}

\theoremstyle{remark}
\newtheorem*{rem*}{Remark}

\newcommand{\sub}{\subseteq}
\newcommand{\R}{\mathbb{R}}
\newcommand{\N}{\mathbb{N}}

\newcommand{\ind}{{\bf 1}}
\DeclareMathOperator{\inter}{int}
\DeclareMathOperator{\dist}{dist}

\DeclareMathOperator{\diam}{diam}

\def\o{\over} 

\title
{Boundary Harnack Principle for fractional powers of Laplacian on
the Sierpi\'nski carpet}
\author{ {\sc Andrzej St\'os}\thanks{
    \emph{2000 MS Classification:} Primary 60J45;
    Secondary 60J35. {\it Key words and phrases}: stable processes,
    Sierpi\'nski carpet, Harnack inequality, Boundary Harnack
    Principle. Research partially supported by KBN grant 2 P03A 041 22 and
    RTN Harmonic Analysis and Related Problems
    contract HPRN-CT-2001-00273-HARP
  } \\
  Laboratoire de Math\'ematiques\\
  Universit\'e Blaise Pascal\\
  24 av. des Landais, 63177 Aubi\`ere Cedex, France \\
  email: stos@math.univ-bpclermont.fr}

\begin{document}
\maketitle

\abstract{We prove the Boundary Harnack Principle related to
  fractional powers of Laplacian for some natural regions in the
  two-dimensional Sierpi\'nski carpet. This is a natual application of
  some more general approach based on the Ikeda-Watanabe formula.}

\begin{center}
  {\bf R\'esum\'e}
\end{center}
Nous pr\'esentons le principe de Harnack \`a la fronti\`ere pour des
puissances fractionaires du laplacien dans les domaines naturels du
tapis de Sierpi\'nski 2-dimensionel. C'est un exemple tr\`es naturel
d'un argument plus g\'en\'eral bas\'e sur la formule d'Ikeda-Watanabe.

\section{Introduction}

Analysis on the Sierpi\'nski carpet (and on a class of similar sets)
has been developing for over ten years (see \cite{BB1}, \cite{BB2} and
references therein). Barlow and Bass showed numerous results including
e.g. the construction of the analogue for the Brownian motion, the
estimates of its transition densities (the heat kernel) and the
Harnack inequality. It is natural to refer to the corresponding
generator as to the Laplacian, even though this is not known whether
this Brownian motion is unique or not.  In this paper we deal with a
fractional power of this Laplacian defined by means of subordination
procedure (see below).  For this operator we give a proof of the
Boundary Harnack 
Principle for some natural regions in the fractal.

In \cite{BSS} (\cite{BSS1}) the Boundary Harnack Principle was
established for {\it cells} in the Sierpi\'nski gasket (or, more
generally, simple nested fractals). The proof in that case resembled
the one for intervals in the real line. In particular, the Boundary
Harnack Principle was a consequence of the (elliptic) Harnack
inequality. This simplification was due to the {\it finite
  ramification} property of the Sierpi\'nski gasket, i.e the fact it
can be disconnected by taking away a finite number of points.  In
particular, the boundary of some natural regions (e.g. small
triangles) is always a  set with a finite number of elements.  
Certainly, the method of \cite{BSS} can
not be carried out to infinitely ramified fractals, such as the
Sierpi\'nski carpet.

In what follows we were influenced by \cite{B} which solves the
problem in the case of Lipschitz domains in $\R^N$.  Our contribution
is a different methodology in proofs which can be described as
follows.  We have no analytic tools and no exact formula for the
Poisson kernel of the ball which are used in \cite{B} (cf.  e.g. Lemma
3 or Lemma 12 in that paper).  Also, a related proof in \cite{SW} uses
theory of smooth functions on $\R^N$.  
Our aim is to present a more general approach
relied on  the Ikeda-Watanabe formula. 
The Sierpi\'nski carpet makes a natural opportunity for application of this
argument. Certainly, the latter depends on the geometric issues,
it seems, however, not to be restricted to this particular fractal.

\section{Preliminaries}
We consider the (unbounded) Sierpi\'nski carpet $F$ which is defined
as follows.  Let $F_0=[0,1]^2$.  Let $A$ be the interior of the {\it
  middle square} of the relative size $1/3$, i.e. $A=(1/3,2/3)^2$.
Set $F_1=F_0\setminus A$.  Then $F_1$ consists of eight closed squares
of side $1/3$.  To obtain $F_2$ we apply subsequently the above
subtraction procedure to these squares in $F_1$, and so on. Set
\begin{displaymath}
  F_\infty= \bigcap_{n=0}^\infty F_n,
  \qquad F=\bigcup_{n=0}^\infty 3^n F_\infty.
\end{displaymath}
We call $F$ the (unbounded) Sierpi\'nski carpet.

By a {\it natural cell} (or simply {\it cell}) we mean the
intersection of $F$ with a square of the form $
[k3^{-n},(k+1)3^{-n}]\times [m3^{-n},(m+1)3^{-n}]$, $k,m,n \in \N$.
The family of cells with sides $3^{-n}$ is denoted by ${\cal S}_n$.

In what follows $D$ always denotes a {\it region} in $F$ i.e.  the
interior of a sum of finite number of natural cells. Since a cell can
be viewed as an union of cells of smaller size, we may and do assume
that $D$ consists of cells which have the same size and disjoint
interiors.  In other words, there exist $n_0,m_0\in \N$, and $S_i\in
{\cal S}_{m_0}$, $i=1,2,...,n_0$ such that
\begin{equation}\label{region}
  D=\inter(\bigcup_{i=1}^{n_0} S_i).
\end{equation}
Note that the interior is taken with respect to the topology of $F$
(inherited from $\R^2$) and since $S_i$ are closed, any two adjacent
cells always make a connected set.  Moreover, the distance between any
two disjoint cells in $D$ is at least $R_1=R_1(D)>0$.  Let
$R_2=3^{-m_0}$ (i. e. $R_2$ is the side of cells in $D$).  Set
$R_0=(1/3)\min(R_1,R_2)$, the number that describes {\it Lipschitz
  character} of $D$.

Notation and conventions. For $x\in F$ and $D \sub F$ we denote
$\delta(x)=\dist(x,\partial D)$.  For $A\sub F$ we write $A^c=F
\setminus A$. By $B(x,r)$ we denote the Euclidean ball (with the
center $x\in F$ and the radius $r>0$) intersected with $F$. For $x,y
\in F$, $|x-y|$ always means the Euclidean distance. Let $d=\dim(F)$
be the Hausdorff dimension of $F$.  By $\mu$ we denote the
$d$-dimensional Hausdorff measure restricted to $F$.  In the sequel
$c$ (without subscripts) denotes a generic constant that depends only
on $F$ and $\alpha$ (see below) and may change its value from one
instance to another.  Constants are numbered consecutively within each
proof.  We write $f(x) \asymp g(x),\;\;x\in F$, to indicate that there
are constants $c_1,c_2 >0$ (independent of $x$) such that $c_1
f(x)\leq g(x) \leq c_2 f(x)$ for all $x\in F$.

To introduce the fractional power of the Laplacian in our framework, we
shortly recall the definition of the {\it $\alpha$-stable process}
from \cite{S} (cf. also \cite{K1},\cite{FJ}).  Let $q(u,x,y)$,
$u>0$, $x,y \in F$, denote transition density (with respect to $\mu$)
of the {\it fractional diffusion} (\cite{Ba}, \cite{BB1}) on $F$.  Set
$\alpha\in(0,2)$ and let $\eta_t(\cdot)$, $t>0$, be a function on
$\R^+$ characterized by its Laplace transform ${\cal
  L}(\eta_t(\cdot))(\lambda) = \exp(-t\lambda^{\alpha/2})$.  (see
\cite{Be} or \cite{BG} for more details and a probabilistic
interpretation).  For $t>0$ and $x,y\in F$ we define
\begin{displaymath}
  p(t,x,y) = \int_0^\infty q(u,x,y)\eta_t(u)du.
\end{displaymath}
By the general theory $p(t,x,y)$ is a transition density of a Markov
process called the subordinate process (see \cite[p. 18]{BG}), which
we denote by $(X_t)_{t>0}$ and call {\it $\alpha$-stable}.
Its generator may be naturally labelled as the $\Delta^{\alpha/2}$.

To simplify the notation, for the rest of the paper we let
$d_\alpha=d+\alpha d_w/2$, where $d_w$ is, in general, a constant
characteristic for the fractal.  For the Sierpi\'nski carpet $d_w
\approx 2.097$.

For a Borel set $B\sub F$ we define {\it exit time} $\tau_B=\inf\{t\ge
0:\;\; X_t\notin B\}$.  Let $u$ be a Borel measurable function $u$ on
$F$, which is bounded from below (above). We say that $u$ is {\it
  $\alpha$-harmonic} in an open set $U\sub F$ if
\begin{displaymath}
  u(x)=E^x u(X(\tau_B)),\qquad x\in B,
\end{displaymath}
for every bounded open set $B$ with the closure $\bar{B}$ contained in
$U$. We say that $u$ is {\it regular $\alpha$-harmonic} in $U$ if
\begin{displaymath}
  u(x)=E^x u(X(\tau_U)),\qquad x\in U.
\end{displaymath}
For a Borel subset $\Omega\sub F$ denote by $\omega^x_\Omega$ the
harmonic measure, i. e. $\omega^x_\Omega(E)=P^x[X_{\tau_\Omega} \in
E]$.

We say that $\Omega \sub F$ has {\it the outer fatness property} (cf.
\cite{BSS}) if there are constants $ c_1=c_1(\Omega)$ and $
r_0=r_0(\Omega)$ such that
\begin{equation}\label{ofp}
  \mu(\Omega^c \cap B(x,r)) \ge c_1 r^d, \qquad x\in \partial
  \Omega,\, r\in (0,r_0). 
\end{equation}
We say that $\Omega$ has {\it the inner fatness property} if there
exist constants $\theta=\theta(\Omega) \in (0,1)$ and
$r_0=r_0(\Omega)$ such that for every $r\in(0,r_0)$ and $Q\in \partial
\Omega$ there is a point $A=A_r(Q) \in \Omega\cap B(Q,r)$ such that
\begin{equation}
  \label{ifp}
  B(A, \theta r)\sub \Omega\cap B(Q,r).
\end{equation}
\begin{rem*}
  Observe that (\ref{ofp}) and (\ref{ifp}) holds for a region $D$.  It
  follows that the Carleson estimate given in Proposition 8.5 of
  \cite{BSS} applies. For the sake of convenience we state it below
  (Lemma \ref{L4}).  Note that if $D$ is a cell of size $3^{-k}$ (or a
  finite union of them) then it satisfies (\ref{ifp}) with
  $r_0=r_0(k)$ and $\theta$ which is an absolute constant, e. g.
  $\theta=1/9$.  We will use this fact without further mention
  dropping from the notation the dependence on $\theta$.
\end{rem*}

\begin{lem}\label{L4}
  Assume $\alpha<2/d_w$.  Let $\Omega \sub F$ be a set satisfying
  (\ref{ifp}).  There exist a constant $c_1=c_1(\theta)$ such that for
  all $Q\in \partial \Omega$ and $r\in(0,r_0/2)$, and functions $u\ge
  0$, regular $\alpha$-harmonic in $\Omega\cap B(Q,2r)$ and satisfying
  $u(x)=0$ on $\Omega^c\cap B(Q,2r)$, we have
  \begin{equation}\label{L4-1}
    u(x)\leq c_1 u(A),\qquad x \in \Omega\cap B(Q,r),
  \end{equation}
  where $A$ is given in (\ref{ifp})
\end{lem}
It can be seen from the proof in \cite{BSS} (cf. also
\cite[(3.29)]{B}) that (\ref{L4-1}) holds for $x\in \Omega\cap
B(Q,5r/4)$, i.e. we have
\begin{equation}\label{L4-2}
  u(x)\leq c_1 u(A),\qquad x \in \Omega\cap B(Q,{5r \o 4}).
\end{equation}
This fact will be invoked later.

Finally, we include the following remark which is due to Prof. Takashi
Kumagai \cite{K2}. The Harnack inequality that we apply here was
proved in \cite{BSS} for $\alpha\in(0,2/d_w)\cup(2d/d_w,2)$. However,
observe that once we have transition density estimates (\cite[Theorem
3.1]{BSS}) then it is relatively easy to deduce the {\it tightness},
i.e. Proposition 4.1 of \cite{CK} for all $\alpha\in(0,2)$. Actually,
this result is contained in \cite[Lemma 4.3]{BSS} (note a different
conventions: $\alpha$ in \cite{CK} means $\alpha d_w/2$ from
\cite{BSS}).  Using this and \cite[Lemma 4.7]{CK} one verifies Lemmas
4.9 - 4.13 of \cite{CK}. Consequently, we can repeat the proof of the
parabolic Harnack inequality \cite[Proposition 4.3]{CK}. This in turn
gives our (elliptic) Harnack inequality for all $\alpha\in (0,2)$.

Unfortunately, in the present paper we have to assume even stronger
restrictions on $\alpha$ (see Lemma \ref{L12}). However, we believe
the restrictions are of the technical nature and once we have the
Harnack inequality for $\alpha\in(0,2)$, the boundary Harnack
Principle holds for the same range of $\alpha$.

\section{Boundary Harnack Principle}

The main result can be stated as follows.
\begin{thm}[Boundary Harnack Principle]\label{bhp}
  Let $\alpha< 2(d-1)/d_w$.  Suppose that $D$ is a region,
  $Q\in\partial D$ and $r\in (0,R_0/2)$.  Then for any functions $u,v
  \ge 0$, positive regular $\alpha$-harmonic in $D\cap B(Q,2r)$ and
  with value $0$ in $D^c\cap B(Q,2r)$, and satisfying
  $u(A_r(Q))=v(A_r(Q))$ we have
  \begin{equation*}
    c_o^{-1} v(x) \le u(x) \le c_o v(x), \qquad x\in D\cap B(Q,r/27),
  \end{equation*}
  where $c_o=c_o(D)$.
\end{thm}

We start the proof by stating some lemmas.  Their assertions have
analogues in \cite{B}. However, there are essential changes in the
argument.  This is required at least for a key step of comparison of
the harmonic measure and the Green function for a region (Lemma
\ref{L12}).  Moreover, the proofs we provide are more elementary in
the sense they rely on basic properties of the process. In particular,
we make use of Ikeda-Watanabe formula and the transition densities
estimates (Proposition 6.1 and Theorem 3.1 in \cite{BSS}). The price
we pay at the moment is the restriction on $\alpha$ (see Lemma
\ref{L12}).

\begin{lem}\label{L10}
  There exist $c_0>0$ such that for any $D$, all $Q \in \partial D$
  and $r\in (0,R_0)$ we have
  \begin{equation*}
    \omega^x_D(B(Q,r)) \ge c_0, \qquad x\in B(Q, r)\cap D.
  \end{equation*}
\end{lem}
\begin{proof}
  Fix $x\in B(Q,r)\cap D$. Recall that $y \to P_D(x,y)$ is the Poisson
  kernel for a region $D$, i.e. the density of $\omega^x(\cdot)$.  By
  \cite[Proposition 6.4]{BSS} and (\ref{ofp}) we get
  \begin{eqnarray*}
    \omega^x_D(B(Q,r)) 
    & \ge &
    P^x[X_{\tau_{B(x,\delta(x)/2)}} \in B(Q,r)\cap D^c] \\
    & \ge & 
    \int_{B(Q,\delta(x))\cap D^c} P_{B(x,\delta(x)/2)} (x,y) d\mu(y) \\
    &\ge & 
    c\delta(x)^{\alpha d_w/2}\int_{B(Q,\delta(x))\cap D^c}
    |x-y|^{-d_\alpha} d\mu(y) \\ 
    &\ge &
    c\delta(x)^{\alpha d_w/2} (2\delta(x))^{-d_\alpha}
    \mu(B(Q,\delta(x)) \cap D^c)\\
    &\ge& c_0,
  \end{eqnarray*}
  which completes the proof.
\end{proof}

Recall that for a region $D$, (\ref{ofp}) and (\ref{ifp}) hold with
some constants $R_0$ and $\theta$.

\begin{lem}\label{L11} 
  Let $\alpha < 2d/d_w$.  There exists a constant $c_1$ such that for
  any region $D$, all $Q\in \partial D$, $r\in (0,R_0)$ and $x\in
  D\setminus B(Q, r)$ we have
  \begin{equation*}
    r^{d-\alpha d_w/2}G_D(x,A_{r/2}(Q)) \le c_1 \omega^x_D(B(Q,r)).
  \end{equation*}
\end{lem}

\begin{proof}
  First we show
  \begin{equation}\label{L11-1}
    \omega^x_D(B(Q,r)) \ge cP^x[T_{B_y}<\tau_D],
  \end{equation}
  where $y=A_{r/2}(Q)$ and $B_y = B(y, {\theta r \o 4})$.  For $x \in
  D$ we have
  \begin{eqnarray*}
    \omega^x_D(B(Q,r)) &\ge& E^x[ \ind_{B(Q,r)} (X_{\tau_D});\,
    T_{B_y}<\tau_D] \\
    & = & E^x[ E^{ X(T_{B_y}) } [\ind_{B(Q,r)} (X_{\tau_D})];\,
    T_{B_y}<\tau_D ] \\
    &\ge& \inf_{w\in B_y} E^w \ind_{B(Q,r)} (X_{\tau_D})
    P^x [T_{B_y}<\tau_D] \\
    &\ge& \inf_{w\in B(Q,r)} \omega^w_D(B(Q,r))P^x[T_{B_y}<\tau_D] \\
    &\ge & c_0 P^x[T_{B_y}<\tau_D],
  \end{eqnarray*}
  where $c_0$ comes from Lemma \ref{L10}.
  
  Now fix $x\in D\setminus B(Q,r)$.  We claim that there exist $c_2$
  such that
  \begin{equation}\label{L11-2}
    c_2G_D(x,y) \delta(y)^{d-\alpha d_w/2} \le
    P^x[T_{B_y}  <\tau_D].
  \end{equation}
  
  To prove our claim observe that $G_D(x,\cdot)$ is $\alpha$-harmonic
  on $D\setminus \{ x\}$ (for $\alpha\not= 2d/d_w$, see e. g.
  \cite{BSS}). Note that $B(y,\delta(y)) \sub B(y, r/2) \sub B(Q,r)$.
  Hence $x\notin B(y,\delta(y))$ and $\overline B(y,\delta(y)) \sub
  D\setminus \{ x\}$.  By the Harnack inequality for the ball
  $B(y,\delta(y))$ we get
  \begin{equation}\label{gd-harnack}
    c_3^{-1}G_D(x,z) \le G_D(x,y) \le c_3 G_D(x,z),\qquad z\in
    B(y,\delta(y)/2).
  \end{equation}
  Since $\theta r/2 <\delta(y)$ we have $B_y\sub B(y,\delta(y)/2)$ and
  hence, by (\ref{gd-harnack}) and the strong Markov property,
  \begin{eqnarray*}
    G_D(x,y)\delta(y)^d &\le& c \theta^{-d} G_D(x,y)\mu(B_y) \\
    &\le& c\int_{B_y}G_D(x,z)d\mu(z) \\
    & = & c G_D\ind_{B_y}(x) \\
    & = & cE^x \left[ \int_0^{\tau_D} \ind_{B_y}(X_s)d s;\,
      T_{B_y}<\tau_D \right] \\
    &=& cE^x[ E^{X(T_{B_y})}\left[
      \int_0^{\tau_D}\ind_{B_y}(X_s)d s \right];\, 
    T_{B_y}<\tau_D ] \\
    & \le & cP^x[ T_{B_y}<\tau_D ] 
    \sup_{ w\in B_y }
    E^w[ \int_0^{\tau_D} \ind_{B_y}  (X_s) d s   ].
  \end{eqnarray*}
  It is easy to see that for $w \in B(y,s)$ we have
  \begin{displaymath}
    \int_{ B(y,s) }{ d\mu(z) \over  |w-z|^{d-\alpha d_w/2}}
    \le \int_{ B(w,2s) }{ d\mu(z) \over  |w-z|^{d-\alpha d_w/2}}
    \le c s^{\alpha d_w/2},
    \qquad s>0,
  \end{displaymath}
  cf. \cite[Lemma 2.1]{BSS}.  It follows that for $w\in B_y$ we have
  \begin{eqnarray*}
    E^w\int_0^{\tau_D} \ind_{B_y}(X_s) d s 
    &\le &
    \int_0^\infty E^w \ind_{B_y}(X_s) d s  \\
    &=& \int_{B_y} \int_0^\infty 
    p(s,w,v) ds \,d\mu(v) \\ 
    &\le& 
    c \int_{B_y} {d\mu(v) \over |v-w|^{d-\alpha d_w/2}} \\
    &\le&
    c\left(\theta \delta(y)\o 4\right)^{\alpha d_w / 2},
  \end{eqnarray*}
  where the last but one inequality is justified by \cite[Lemma
  5.3]{BSS}. Note that this is the only place where we used
  $\alpha<d_s$.  The claim follows.
  
  Since $\theta r/2 \le \delta(y) \le r/2$ (i.e. $\delta(y) \asymp
  r$), (\ref{L11-1}) and (\ref{L11-2}) imply the assertion of the
  lemma.
\end{proof}

\begin{lem}\label{L12}
  If $\alpha < 2(d-1)/d_w$ then there exists a constant $c_1$ such
  that for any $D$, all $Q\in \partial D$ and $r\in (0,R_0/2)$ we have
  \begin{equation*}
    \omega^x_D(B(Q,r)) \le c_1 r^{d-\alpha d_w/2} G_D(x,A_{r /2}(Q)),
    \qquad x\in D\setminus B(Q, 2r).
  \end{equation*}
\end{lem}

\begin{proof}
  Fix $x\in D\setminus B(Q,2r)$.  It can be observed that the harmonic
  measure does not charge $\partial D$.  Indeed, it is enough to adapt
  Lemma 6 of \cite{B} with {\it outer cone property} replaced by
  (\ref{ofp}).  For the sake of reader's convenience we sketch the
  argument.  Denote $\tau_x = \tau_{B(x,\delta(x)/3)}$. Then, by the
  strong Markov property,
  \begin{equation*}
    \omega^x_D(\partial D) = P^x[X_{\tau_x} \in \partial D] 
    +E^x[\omega^{X_{\tau_x}}_D;\, X_{\tau_x}\in D] =: p_0(x)+r_0(x).
  \end{equation*}
  Define inductively
  \begin{equation*}
    p_{k+1}(x) = E^x[p_k(X_{\tau_x});\, X_{\tau_x}\in D],
  \end{equation*}
  \begin{equation*}
    r_{k+1}(x) = E^x[r_k(X_{\tau_x});\, X_{\tau_x}\in D].
  \end{equation*}
  Then $r_k=p_{k+1}+r_{k+1}$, $k=0,1,..$, and
  \begin{equation}
    \label{L6-1}
    \omega^x_D(\partial D) = p_0(x)+p_1(x)+...+p_k(x)+r_k(x),\quad x\in
    D,\, k=0,1,...
  \end{equation}
  Let $x_0\in \partial D$ be such that $|x_0-x|=\delta(x)$.  By
  \cite[Proposition 6.4]{BSS} and (\ref{ofp}) we get
  \begin{eqnarray*}
    P^x[X_{\tau_x} \in D^c] 
    &\ge& 
    P^x[X_{\tau_x} \in B(x_0,\delta(x))\cap D^c] 
    \\& \ge &
    c\delta(x)^{\alpha d_w/2} \int_{B(x_0,\delta(x))\cap D^c} 
    {d\mu(y) \o |x-y|^{d_\alpha}} 
    \\&\ge& 
    {c\delta(x)^{\alpha d_w/2} \o (2\delta(x))^{d_\alpha}}
    \mu(B(x_0,\delta(x))\cap D^c)
    \\ &\ge & c_0,
  \end{eqnarray*}
  for each $x\in D$.  Consequently,
  \begin{equation*}
    \sup_{x\in D}r_{k+1}(x) \le (1-c_0)\sup_{x\in D}r_k(x) \le
    (1-c_0)^{k+1} \longrightarrow 0,\quad k\to \infty.
  \end{equation*}
  From (\ref{L6-1}) it follows that
  \begin{equation*}
    \omega^x_D(\partial D) = \sum_{k=0}^\infty p_k(x).
  \end{equation*}
  Since $\mu$ does not charge $\partial D$ we immediately get $p_k(x)
  = 0$, $x\in D$, $k=0,1,..$ (see also the remark after Corollary 6.2
  in \cite{BSS}). This gives our claim.
  
  Now, since $\omega^x_D(\partial D) =0$, from the Ikeda-Watanabe
  formula (see also \cite[(51)]{BSS}) we have
  \begin{eqnarray*}
    \omega^x_D(B(Q,r)) 
    & = & \int_{B(Q,r)\cap D^c} P_D(x,y) d\mu(y)\\
    & \asymp &
    \int_{B(Q,r)\cap D^c} \int_D
    {G_D(x,z) \over |z-y|^{d_\alpha}}d\mu(z)d\mu(y) \\
    &=& 
    \left(\int_{D\setminus B(Q,{5r\o 4})} +  \int_{D\cap B(Q,{5r\o 4})} \right)
    \left[
      G_D(x,z) 
      \int_{B(Q,r)\cap D^c} { d\mu(y) \over |z-y|^{d_\alpha}} 
    \right] 
    d\mu(z)\\
    &=& J_1+J_2.
  \end{eqnarray*}
  First we deal with the integral $J_1$. Let $A_0=A_{r/2}(Q)$.  Then
  we have $|z-y| \ge r/4 $ and so $|z-A_0| \le |z-y|+|y-A_0| \le
  |z-y|+(3/2)r \le |z-y|+6|z-y| = 7|z-y|$.  It follows that
  \begin{equation*}
    \int_{B(Q,r)\cap D^c} {d\mu(y) \over |z-y|^{d_\alpha}} 
    \le { c\over |z-A_0|^{d_\alpha} } \mu(B(Q,r)) \asymp 
    { cr^d \over |z-A_0|^{d_\alpha} }
  \end{equation*}
  and
  \begin{equation}\label{I2}
    J_1 \le cr^d
    \int_{D\setminus B(Q,5r/4)} {G_D(x,z) \o |z-A_0|^{d_\alpha}}
    d\mu(z).
  \end{equation}
  Denote $B_0=B(A_0, \theta r/2)$.  For the Poisson kernel of the ball
  $B_0$ by \cite[Proposition 6.4]{BSS} we have
  \begin{equation*}
    P_{B_0}(A_0,z) \ge c { (\theta r/2)^{\alpha d_w/2} \o
      |z-A_0|^{d_\alpha} }, \quad z\in B_0^c. 
  \end{equation*} 
  By rearranging and putting this into (\ref{I2}) we obtain
  \begin{equation*}
    J_1\le cr^{d-\alpha d_w/2} \int_{B_0^c} P_{B_0}(A_0,z)G_D(x,z)d\mu(z).
  \end{equation*}
  Since $z\to G_D(x,z)$ is regular $\alpha$-harmonic on $B_0$, the
  last integral does not exceed $G_D(x,A_0)$.  Remark that the
  integral is not necessarily equal to $G_D(x,A_0)$, since we do not
  know whether the process hits the boundary of $B_0$; however, we do
  not need this fact and the equality.  Finally,
  \begin{equation}\label{j1}
    J_1 \le cr^{d-\alpha d_w/2} G_D(x,A_{r/2}(Q)),
  \end{equation}
  as desired.
  
  To deal with the integral $J_2$ observe that
  \begin{equation*}
    \int_{B(Q,r)\cap D^c} {d\mu(y) \o |z-y|^{d_\alpha}}
    \le 
    \int_{B(z,\delta(z))^c} {d\mu(y) \o |z-y|^{d_\alpha}}
    \le 
    c\delta(z)^{-\alpha d_w/2},
  \end{equation*}
  where the last inequality is justified by Lemma 2.1 of \cite{BSS}.
  Since $z\mapsto G_D(x,z)$ is regular $\alpha$-harmonic on $D\cap
  B(Q,2r)$, from (\ref{L4-2}) it follows that
  \begin{eqnarray}
    \nonumber J_2 &\le& 
    c\int_{D\cap B(Q, 5r/4)} G_D(x,z)\delta(z)^{-\alpha d_w/2} d\mu(z)\\
    &\le& 
    cG_D(x,A_{5r/4}(Q))\int_{D\cap B(Q, 5r/4)} \delta(z)^{-\alpha
      d_w/2}d\mu(z). \label{j2a}
  \end{eqnarray}
  
  We have $|A_{5r/4}-A_0| \le |A_{5r/4}-Q|+|Q-A_0| \le 5r/4 +r/2 \le c
  (\theta r/2)$.  By \cite[Lemma 7.6]{BSS} with $x_1=A_{5r/4}$ and
  $x_2=A_0=A_{r/2}(Q)$ we obtain
  \begin{equation}\label{j2b}
    G_D(x,A_{5r/4}(Q)) \le c G_D(x,A_{r/2}(Q)).
  \end{equation}
  Now, it is enough to estimate
  \begin{equation*}
    \int_{D\cap B(Q,5r/4)} \delta(z)^{-\alpha d_w/2}d\mu(z).
  \end{equation*}
  
  Let $k_o\in \N$ be such that $3^{-k_o-1} < 5r/4 \le 3^{-k_o}$. Then,
  clearly, $r \asymp 3^{-k_o}$. Let $H_0$ be the union of cells $S$
  that satisfy
  \begin {itemize}
  \item[(a)] $S\in {\cal S}_{k_o}$,
  \item[(b)] $S \sub \overline{D}$,
  \item[(c)] $\partial S \cap \partial D \not= \emptyset$,
  \item[(d)] $S\cap B(Q,5r/4) \not= \emptyset$.
  \end{itemize}
  In other words $H_0$ is a covering of $D\cap B(Q,5r/4)$ by {\it
    smallest} cells adjacent to $\partial D$. Define $H_k$,
  $k=1,2,...$, in the same way as $H_0$ but with (a) replaced by
  $S\in{\cal S}_{k_o+k}$ and (d) replaced by $S\sub H_0$.  Thus, $H_k$
  is a layer of cells of side $3^{-k-k_o}$ adjacent to $\partial D
  \cap \partial H_0$.  Then, there is at most $h_k=2.3^k+1$ cells in
  $H_k$, $k=1,2,...$ (this may happen when $H_0$ consists of three
  cells, i.e. $Q\in \delta D$ is a corner point).  Let
  $R_k=H_k\setminus H_{k+1}$. Then $z\in R_k$ implies $\delta(z)\ge
  3^{-(k_o+k+1)} \ge cr3^{-k}$. It follows that
  \begin{eqnarray}
    \label{j2c} \int_{D\cap B(Q,5r/4)} \delta(z)^{-\alpha d_w/2}d\mu(z)
    & \le &
    \sum_{k=0}^\infty \int_{R_k}\delta(z)^{-\alpha d_w/2}d\mu(z) 
    \\ \nonumber & \le  &
    c\sum_{k=0}^\infty  (3^{-k} r)^{-\alpha d_w/2}
    \mu(R_k) 
    \\  \nonumber & \le & 
    c r^{-\alpha d_w/2}\sum_{k=0}^\infty
    3^{k \alpha d_w/ 2} (3^{-k} r)^d h_k 
    \\ \nonumber & \le & 
    c r^{d-\alpha d_w/2}\sum_{k=0}^\infty
    3^{k(\alpha d_w/2-d+1)} 
    \\ \nonumber & \le & 
    cr^{d-\alpha d_w/2}, 
  \end{eqnarray}
  provided $\alpha <2(d-1)/d_w$.  Combining (\ref{j1}), (\ref{j2a}),
  (\ref{j2b}) and (\ref{j2c}) we get the assertion.
\end{proof}

\begin{rem*}
  In our particular case $2(d-1)/d_w\approx 0.851$.
\end{rem*}

\begin{proof}[Proof of Theorem \ref{bhp}]
  This is based on a general idea of the proof of Lemma 13 from
  \cite{B}.  Since the context is different, we present a version
  adapted to our needs. The argument goes the following way. First, we
  introduce the basic geometrical objects and notations.  Then, the
  first step of the proof is to establish the comparability of the
  harmonic measures of the region $\Delta$ and of its propper subset
  $B_1$ (see below). This is given in (\ref{bhp1}) which is a key
  ingredient in the proof.  Then we decompose the functions to be
  compared into two parts (\ref{split}). In Steps 2 and 3 we prove the
  inequality for each of these parts: (\ref{bhp10}) and (\ref{step3})
  respectively.  Step 2 is the crucial one and it uses (\ref{bhp1});
  Step 3 is covered by the Poisson kernel estimates and the (usual)
  Harnack inequality.
  
  \begin{center}
    \begin{picture}(283,283)(-40,-40)
      \linethickness{0.2pt}
      \multiput(0,0)(3,0){82}{\line(0,1){243}}
      \multiput(0,0)(0,3){82}{\line(1,0){243}}      
      \put(81,81){\rule{81pt}{81pt}}
      \multiput(0,0)(0,27){9}{\multiput(9,9)(27,0){9}{\rule{9pt}{9pt}}}
      \multiput(0,0)(0,81){3}{\multiput(27,27)(81,0){3}{\rule{27pt}{27pt}}}
      \multiput(0,0)(0,9){27}{\multiput(3,3)(9,0){27}{\rule{3pt}{3pt}}}
      
      \put(37,-6){$_Q$}
      \put(40.4,0){\circle*{2}}
      \put(250,121){$3^{-N}$}
      
      \put(-30,-30){\line(1,0){192}}
      \put(162,-29.9){\vector(0,1){30}}
      \put(-30,-30){\line(0,1){192}}
      \put(-30,162){\vector(1,0){30}}
      \put(133,-26){$_\Omega$}
      \put(-27,134){$_\Omega$}
      
      \put(-15,159){\vector(1,0){15}}
      \put(-15,153){$_{A_i}$}
      \put(159,-15){\vector(0,1){15}}
      \put(149,-13){$_{A_i}$}
      
      \put(-20,-20){\line(1,0){101}}
      \put(81,-20){\vector(0,1){20}}
      \put(-20,-20){\line(0,1){47}}
      \put(-20,27){\vector(1,0){20}}
      \put(-18,20){$_{\Omega_2}$}

      \put(-30,162){\line(0,1){9}}
      \put(-30,171){\vector(1,0){30}}
      \put(-29,164){\tiny${\tilde B_i}$}
      
      \put(162,-30){\line(1,0){9}}
      \put(171,-30){\vector(0,1){30}}
      \put(162.5,-24){\tiny ${\tilde B_i}$}
      
      \put(234,-10){\line(1,0){9}}
      \put(243,-12){\line(0,1){4}}
      \put(234,-12){\line(0,1){4}}
      \put(238,-17){{\tiny $\tilde r$}}
    \end{picture}
  \end{center}
  
  Let $N\in \N$ be such that $3^{-N} \le r<3^{-N+1}$.  For $Q\in
  \partial D$ let $S^\nu_i(Q)$, $i=1,2$, be cells from ${\cal
    S}_{N+i}$ such that $Q\in S^\nu_i(Q)\sub D$.  There can be one,
  two or three such cells indexed by $\nu$.  Define
  \begin{equation*}
    \Omega_i=\inter \left(\bigcup_\nu S^\nu_i(Q)\right), \qquad i=1,2.
  \end{equation*}
  If the union above consists of the single $S^1_i(Q)$ then we set
  \begin{equation*}
    \Omega_i= \inter \left(S^1_i(Q)\cup\bigcup_{\nu=1}^2
      N_i^\nu(Q)\right), \qquad i=1,2. 
  \end{equation*}
  where $N_i^\nu$ are the {\it neighbours} of $S_i^1(Q)$, i.e.  cells
  satisfying
  \begin{itemize}
  \item[(i)] $N_i^\nu\in {\cal S}_{N+i}$ and $N_i^\nu \sub \overline
    D$,
  \item[(ii)] $\partial N_i^\nu \cap \partial D \cap \partial S^1_i(Q)
    \not= \emptyset$ (recall that cells are closed).
  \end{itemize}
  Finally, denote $\Omega=\Omega_1$.
  
  Set $\tilde r=3^{-N-3}$ and let $A \in \Omega$ be a point such that
  $\dist(A,D^c) = 3\tilde r$ and $\dist(A,\Omega_2)=\tilde r$
  (clearly, $A$ is not unique).

  \begin{rem*}
    In the course of the proof it is convenient to identify $A$ with
    $A_r(Q)$ from the hypothesis of our theorem. Note that there is no
    loss of generality; indeed, by \cite[Lemma 7.6]{BSS} we have
    $u(A_r(Q)) \asymp u(\tilde A_r(Q))$ for any harmonic function $u$
    satisfying hypothesis of Theorem 3.1 and points $A_r(Q)$, $\tilde
    A_r(Q)$ of the inner fatness property. Actually, this is the
    reason we can use our our definition of $A$ and $A_r(Q)$ without
    determining uniquely the points.
  \end{rem*}
  
  Let $\tilde B_i \in {\cal S}_{N+3}$, $i=1,2,..,n_0(\Omega)$, are
  cells satisfying $\tilde B_i\sub \overline{D}\cap \Omega^c$ and
  $\partial \tilde B_i \cap \partial\Omega\not= \emptyset$.  Since
  $18\le n_0(\Omega) \le 54$, we drop the dependence $n_0$ on $\Omega$
  without further mention.  Set $\tilde B_1$ to be one of $\tilde B_i$
  satisfying additionally $\dist(\tilde B_1,\partial D)\ge 8\tilde r$.
  Let $S_i$ be the mid-point of the line segment
  $\partial\Omega\cap\partial \tilde B_i$; if the set consists of one
  point $\{ x_o \}$ then let $S_i=x_o$( a {\it vertex point}).  Let
  $B_i=B(S_i,\tilde r\sqrt{2})$ and
  \begin{equation*}
    \Delta= \bigcup_i B_i\cap D \cap \Omega^c. 
  \end{equation*}
  Let $A_i\in \Omega$, $i=1,2,..n_0$, be the point such that
  $|A_i-S_i|= \dist(A_i, \delta(\Omega))=\tilde r/3$, provided $S_i$
  is not a vertex point of $\Omega$, and $|A_i-S_i|= \tilde
  r\sqrt{2}/3$ in the opposite case.  $\dist(A_i,
  \delta(\Omega))=\tilde r/3$.  Since $\dist(\tilde B_1,\partial D)\ge
  8\tilde r$ then there exists a cell, denoted by $T$, such that $T\in
  {\cal S}_{N+4}$, $T\sub D\setminus (\Omega \cup \Delta)$,
  $\dist(T,D^c)\ge 8 \tilde r$ and $\dist(T,B_1)\leq \tilde r$.
  
  {\sc Step 1.}  Let $\theta=1/9$. Then if $x\in B(A_i, \theta \tilde
  r\sqrt{2}/2)$ then $|x-S_i| \le |x-A_i|+|A_i-S_i|\le \tilde r
  \sqrt{2} /18 +\tilde r\sqrt{2}/3 \le \tilde r\sqrt{2}/2$, which
  yields $B(A_i, \theta \tilde r\sqrt{2}/2) \sub \Omega\cap B(S_i,
  \tilde r\sqrt{2}/2)$. In other words, $A_i$ can be regarded as
  $A_{\tilde r \sqrt 2/2}(S_i)$ in the inner fatness property
  (\ref{ifp}) for $\Omega$.  It follows that by Lemmas \ref{L11} and
  \ref{L12} applied to $\Omega$ and $B_i$ we get
  \begin{equation*}
    (\tilde r \sqrt{2})^{d-\alpha d_w/2} G_\Omega(z,A_i) \asymp
    \omega^z_\Omega(B_i), 
    \qquad z\in \Omega\setminus B(S_i,2\tilde r\sqrt{2}).
  \end{equation*}
  
  For the rest of the proof fix $x\in \Omega_2$. Then $|x-S_i| \ge
  6\tilde r$, $i=1,2,...,n_0$, and hence
  \begin{equation*}
    \tilde r^{d-\alpha d_w/2} G_\Omega(x,A_i) \asymp \omega^x_\Omega(B_i).
  \end{equation*}
  
  Recall $\dist(A,D^c) = 3\tilde r$.  Since $\dist(A_i,\partial\Omega)
  = \tilde r/3$, $|A_i-A|\le \diam(\Omega)\le c(\tilde r/3)$ and
  $G_\Omega(x,\cdot)$ is regular $\alpha$-harmonic in $B(A_i,\tilde
  r/3) \cup B(A,\tilde r/3)$, by Harnack inequality (\cite[Lemma
  7.6]{BSS}) we obtain
  \begin{equation}\label{bhp2}
    G_\Omega(x,A_i) \asymp G_\Omega(x,A).
  \end{equation}
  It follows that
  \begin{equation}\label{bhp1}
    \begin{array}{ccl}
      \omega^x_\Omega(\Delta)
      &\le&
      \sum_{i=1}^{n_0} \omega^x_\Omega(B_i) \\
      &\asymp&
      \tilde r^{d-\alpha d_w/2}\sum_{i=0}^{n_0} G_\Omega(x,A_i) \\
      & \asymp & \tilde r^{d-\alpha d_w/2} G_\Omega(x,A_1) \\
      & \asymp & \omega^x_\Omega(B_1).
    \end{array}
  \end{equation}
  
  {\sc Step 2.} Let $u_1,\, u_2$ be functions such that
  \begin{equation}\label{split}
    \begin{array}{ccccc}
      u_1(y)= \left\{
        \begin{array}{cc}
          u(y),&  y \in \Delta, \\
          0,   &  y \in \Omega^c\setminus \Delta,
        \end{array}
      \right.
      &  & & &
      u_2(y)= \left\{
        \begin{array}{cc}
          0,&  y \in \Delta \\
          u(y),   &  y \in \Omega^c\setminus \Delta,
        \end{array}
      \right.
    \end{array}
  \end{equation}
  and $u_1$ and $u_2$ are regular $\alpha$-harmonic in $\Omega$.  Note
  that $u_1, u_2 \ge 0$ and $u_1+u_2=u$.  Analogously we define $v_1$
  and $v_2$.
  
  By (\ref{L4-1}) and (\ref{bhp1}) we obtain
  \begin{equation}\label{bhp3}
    \begin{array}{ccl}
      u_1(x) 
      &=& 
      E^x[u(X_{\tau_\Omega});\;X_{\tau_\Omega}\in \Delta] 
      \\ & \le &
      \sup \{ u(z);\; z\in \Delta \} \omega^x_\Omega(\Delta)
      \\ & \le &
      cu(A)\omega^x_\Omega(\Delta) 
      \\ &\le&  
      cu(A)\omega^x_\Omega(B_1).
    \end{array}
  \end{equation}
  
  Since $\dist(A\cup B_1,\partial D) \ge \tilde r$ and for $y\in B_1$
  we have $\dist(A,y)\le \diam(\Omega)+\diam(B_1) \le c\tilde r$, from
  \cite[Lemma 7.6]{BSS} it follows that
  \begin{equation*}
    v_1(y)=v(y) \ge cv(A), \qquad y\in B_1.
  \end{equation*} 
  Consequently, we have
  \begin{eqnarray*}
    v_1(x) &=& E^x[v(X_{\tau_\Omega});\; X_{\tau_\Omega}\in \Delta ] \\
    &\ge & 
    E^x[v(X_{\tau_\Omega});\; X_{\tau_\Omega}\in B_1 ]\\
    &\ge&
    cv(A)\omega^x_\Omega(B_1).
  \end{eqnarray*}
  Combinig this and (\ref{bhp3}) we get
  \begin{equation}\label{bhp10}
    u_1(x) \le cv_1(x) \le cv(x).
  \end{equation}
  
  {\sc Step 3.} Now, let $K=\Omega\cup \Delta\cup (D^c\cap B(Q,2r))$.
  Clearly, $\bigcup_i \tilde B_i \sub \Delta$. So if $z\in
  D\setminus(\Omega\cup \Delta)$ then $\dist(z,\Omega) \ge \tilde r$.
  Hence, for $z\in \Omega$ and $y\in K^c$ we have $|y-z| \asymp
  |y-Q|$.  Therefore, by the Ikeda-Watanabe formula
  \begin{eqnarray*}
    u_2(x) & = & \int_{K^c}P_\Omega(x,y)u(y)d\mu(y)\\
    &\asymp & \int_{K^c} \left( \int_\Omega
      G_\Omega(x,z)|z-y|^{-d_\alpha} d\mu(z)\right)u(y)d\mu(y)\\
    &\asymp&
    \int_{K^c} \left( \int_\Omega
      G_\Omega(x,z) d\mu(z)\right) u(y)|y-Q|^{-d_\alpha} d\mu(y)\\
    &=&
    E^x \tau_\Omega \int_{K^c}
    u(y)|y-Q|^{-d_\alpha} d\mu(y)
  \end{eqnarray*}
  From this and the analogous relation for $v_2$ it follows that
  \begin{equation}\label{bhp8}
    u_2(x)/u_2(A) \asymp E^x\tau_\Omega/E^{A}\tau_\Omega
    \asymp v_2(x)/v_2(A).
  \end{equation}  
  We claim that
  \begin{equation}\label{bhp9}
    v_2(A) \ge cv(A).
  \end{equation} 
  Indeed, recall that $T\cap\Delta=\emptyset$ and we have
  \begin{equation}\label{bhp6}
    v_2(A) \ge  
    E^A[ v(X_{\tau_\Omega});\; X_{\tau_\Omega}\in T] 
    \ge 
    \inf_{z\in T}v(z)\omega^A_\Omega(T).
  \end{equation}
  Since $\dist(A\cup T,\partial D)\ge 3\tilde r$ and $\dist(A,T)\le
  c\tilde r$, by the Harnack inequality we have
  \begin{equation}\label{bhp7}
    v(z) \asymp v(A), \qquad z\in T.
  \end{equation}
  Moreover, $\diam(\Omega)\asymp\diam(T)\asymp\dist(\Omega,T)\asymp
  \tilde r$ yields $|y-z| \asymp \tilde r$, $y\in \Omega$, $z\in T$.
  Hence, by \cite[Proposition 4.4]{BSS}
  \begin{eqnarray*}
    \omega^A_\Omega(T) 
    &\asymp&
    \int_T \int_\Omega {G_\Omega(A,y) \o |y-z|^{d_\alpha}} d\mu(y)d\mu(z)  \\
    &\asymp&
    \tilde r^{-d_\alpha}\int_T \int_\Omega G_\Omega(A,y) d\mu(y)d\mu(z) \\
    &=& \mu(T)  \tilde r^{-d_\alpha} E^A\tau_\Omega \\
    &\ge & c\tilde r^{-\alpha d_w/2}E^A\tau_{B(A,\tilde r)} = c_1,
  \end{eqnarray*}
  where $c_1$ is independent of $\Omega$, $T$, $r$, etc.
  Putting this and (\ref{bhp7}) into (\ref{bhp6}) we get our claim.
  
  Denote the last quotient in (\ref{bhp8}) by $q_o$.  Then, by
  (\ref{bhp8}), definition of $u_2$, the assumption $u(A)=v(A)$ and
  (\ref{bhp9}),
  \begin{eqnarray}\label{step3}
    &  u_2(x) \le cq_o u_2(A) \le cq_o u(A) = cq_ov(A) & \\
    & \le c q_o v_2(A) =c v_2(x) & x\in \Omega_2.\nonumber
  \end{eqnarray}
  Together with (\ref{bhp10}) and the symmetry this ends the proof.
\end{proof}
\begin{rem*}
  Although the proof relies on particular geometric properties of the
  Sierpi\'nski carpet, we believe that this argument can be carried
  out to a slightly wider context, e.g. to generalized Sierpi\'nski
  carpets.
\end{rem*}

\subsection*{Acknowledgements}    
Part of this work was done during post-doctoral fellowship in
Universit\'e d'Angers. Financial support by the region region Pays de
la Loire is gratefully acknowledged.

\end{document}